\documentstyle[12pt]{amsart}
\begin{document}
  \title{Harmonic maps and biharmonic Riemannian submersions}
  \title[Harmonic maps and biharmonic Riemannian submersions]
  {Harmonic maps and biharmonic Riemannian submersions}
   \author{Hajime Urakawa}
  \address{Global Learning Center, 
  Tohoku University, Kawauchi 41, Sendai 980-8576, Japan}
  \curraddr{Graduate School of Information Sciences, Tohoku University, Aoba 6-3-09, Sendai 989-8579, Japan}
  \email{urakawa@@math.is.tohoku.ac.jp}
    \keywords{Riemannian submersions, 
    harmonic map, biharmonic map}
  \subjclass[2000]{primary 58E20, secondary 53C43}
  \thanks{
  Supported by the Grant-in-Aid for the Scientific Reserch, (C) No. 25400154, Japan Society for the Promotion of Science. 
  }
\maketitle
\begin{abstract} 
Characterizations for Riemannian submersions to be harmonic or biharmonic are shown. Examples of 
  biharmonic but not harmonic Riemannian submersions are shown.
    \end{abstract}
\numberwithin{equation}{section}
\theoremstyle{plain}
\newtheorem{df}{Definition}[section]
\newtheorem{th}[df]{Theorem}
\newtheorem{prop}[df]{Proposition}
\newtheorem{lem}[df]{Lemma}
\newtheorem{cor}[df]{Corollary}
\newtheorem{rem}[df]{Remark}
\section{Introduction}
Variational problems play central roles in geometry;\,Harmonic map is one of important variational problems which is a critical point of the energy functional 
$E(\varphi)=\frac12\int_M\vert d\varphi\vert^2\,v_g$ 
for smooth maps $\varphi$ of $(M,g)$ into $(N,h)$. The Euler-Lagrange equations are given by the vanishing of the tension filed 
$\tau(\varphi)$. 
In 1983, J. Eells and L. Lemaire \cite{EL1} extended the notion of harmonic map to  
biharmonic map, which are, 
by definition, 
critical points of the bienergy functional
\begin{equation}
E_2(\varphi)=\frac12\int_M
\vert\tau(\varphi)\vert^2\,v_g.
\end{equation}
After G.Y. Jiang \cite{J} studied the first and second variation formulas of $E_2$, 
extensive studies in this area have been done
(for instance, see 
\cite{CMP}, \cite{LO2},  \cite{MO1}, \cite{OT2}, \cite{S1},
\cite{IIU2}, \cite{IIU},  \cite{II}, 
 etc.). Notice that harmonic maps are always biharmonic by definition. 
B.Y. Chen raised (\cite{C}) so called B.Y. Chen's conjecture and later, R. Caddeo, S. Montaldo, P. Piu and C. Oniciuc raised (\cite{CMP}) the generalized B.Y. Chen's conjecture. 
\par
\textbf{B.Y. Chen's conjecture:}\par
 {\em Every biharmonic submanifold of the Euclidean space ${\mathbb R}^n$ must be harmonic (minimal).}
\par
\textbf{The generalized B.Y. Chen's conjecture:}\par
{\em Every biharmonic submanifold of a Riemannian manifold of non-positive curvature 
must be harmonic (minimal).}
\vskip0.6cm\par
For the generalized Chen's conjecture, 
Ou and Tang gave (\cite{OT}, \cite{OT2}) a counter example in a Riemannian manifold of negative curvature. 
For the Chen's conjecture, affirmative answers were known for the case of 
surfaces in the three dimensional Euclidean space (\cite{C}), 
and the case of 
hypersurfaces of the four dimensional Euclidean space (\cite{HV}, \cite{D}). 
K. Akutagawa and S. Maeta gave (\cite{AM}) showed
a supporting evidence to the Chen's conjecture: 
{\em
Any complete regular biharmonic submanifold of the Euclidean space 
${\mathbb R}^n$ is harmonic (minimal). 
}
The affirmative answers to the generalized B.Y. Chen's conjecture were shown 
(\cite{NU1}, \cite{NU2}, \cite{NUG}) 
under the $L^2$-condition and completeness of $(M,g)$. 
\vskip0.1cm\par
In \cite{U1}, we treated with a principal $G$-bundle over a Riemannian manifold, and 
showed the following two theorems:  
\vskip0.3cm\par
{\bf Theorem 1}
{\em Let $\pi:\,(P,g)\rightarrow (M,h)$ be a principal 
 $G$-bundle 
over a Riemannian manifold $(M,h)$ 
with non-positive Ricci curvature. 
 Assume $P$ is compact so that $M$ is also compact. 
If the projection $\pi$ is biharmonic, then it is harmonic. }
\vskip0.3cm\par
{\bf Theorem 2} 
{\em 
Let $\pi:\,(P,g)\rightarrow (M,h)$ be a principal $G$-bundle over a Riemannian manifold with 
non-positive Ricci curvature. 
Assume that $(P,g)$ is a non-compact complete Riemannian manifold, and 
the projection $\pi$ has both finite energy $E(\pi)<\infty$ and finite bienergy $E_2(\pi)<\infty$. 
If $\pi$ is biharmonic, then it is harmonic. 
}
\vskip0.3cm\par
We give two comments on the above theorems: For the generalized B.Y. Chen's conjecture, non-positivity of the sectional curvature of the ambient space 
of biharmonic submanifolds is necessary. However, it should be emphasized that 
for the principal $G$-bundles, we need not the assumption of non-positivity of the sectional curvature. We only assume 
{\em non-positivity of the Ricci curvature} of the domain manifolds in the proofs of Theorems 1 and 2. Second, in Theorem 2, finiteness of the energy and bienergy is necessary. Otherwise, one can see the following counter examples due to
Loubeau and Ou (\cite{LOu}): 
\vskip0.3cm\par
{\bf Example 1} (cf. \cite{BK}, \cite{LOu}, p. 62) 
The inversion in the unit sphere
$\phi:\,{\mathbb R}^n\backslash 
\{o\}\ni {\bf x}\mapsto 
\frac{\bf x}{\vert {\bf x}\vert^2}\in {\mathbb R}^n$ is biharmonic  if $n=4$. It is not harmonic 
since $\tau(\phi)=-\frac{4\,{\bf x}}{\vert{\bf x}\vert^4}$.
\vskip0.3cm\par 
{\bf Example 2} (cf. \cite{LOu}, p. 70) 
Let $(M^2,h)$ be a Riemannian surface, 
and let 
$\beta:\,M^2\times {\mathbb R}\rightarrow {\mathbb R}^{\ast}$ and $\lambda:\,{\mathbb R}\rightarrow {\mathbb R}^{\ast}$ 
be two positive $C^{\infty}$ functions. 
Consider the projection 
$\pi:\,(M^2\times {\mathbb R}^{\ast}, 
g=\lambda^{-2}\,h+\beta^2\,dt^2)\ni (p,t)\mapsto p\in (M^2,h)$. 
Here, we take $\beta= c_2\,e^{\int f(x)\,dx}$, 
$f(x)=\frac{-c_1\,(1+e^{c_1x})}{1-e^{c_1x}}$ with 
$c_1,\,c_2\in {\mathbb R}^{\ast}$, and 
$(M^2,h)=({\mathbb R}^2,dx^2+dy^2)$. Then, 
\vskip0.2cm\par\noindent
$
\pi:\,({\mathbb R}^2\times {\mathbb R}^{\ast},dx^2+dy^2+\beta^2(x)\,dt^2)\ni (x,y,t)\mapsto (x,y)\in ({\mathbb R}^2,dx^2+dy^2)
$
\vskip0.2cm\par\noindent
gives a family of {\em proper biharmonic} 
(i.e., biharmonic but not harmonic)
Riemannian submersions.   
\vskip0.3cm\par
In this paper, we will treat with a more general setting of 
Riemannian submersion 
$\pi:\,\,(P, g)\rightarrow (M,h)$ 
with a $S^1$ fiber over a compact Riemannian manifold 
$(M,h)$. 
We first derive the tension field $\tau(\pi)$ and 
the bitension field $\tau_(\pi)$. As a corollary of our main theorem, we show a characterization theorem for a Riemannian submersion 
$\pi:\,(P,g)\rightarrow (M,h)$ 
over a compact K\"{a}hler-Einstein manifold $(M,h)$, 
 to be {\em biharmonic}.  
\vskip0.6cm\par
\section{Preliminaries}
\subsection{Harmonic maps and 
biharmonic maps.}
We first prepare the materials for the first and second variational formulas for the bienergy functional and biharmonic maps. 
Let us recall the definition of a harmonic map $\varphi:\,(M,g)\rightarrow (N,h)$, of a compact Riemannian manifold $(M,g)$ into another Riemannian manifold $(N,h)$, 
which is an extremal 
of the {\em energy functional} defined by 
$$
E(\varphi)=\int_Me(\varphi)\,v_g, 
$$
where $e(\varphi):=\frac12\vert d\varphi\vert^2$ is called the energy density 
of $\varphi$.  
That is, for any variation $\{\varphi_t\}$ of $\varphi$ with 
$\varphi_0=\varphi$, 
\begin{equation}
\frac{d}{dt}\bigg\vert_{t=0}E(\varphi_t)=-\int_Mh(\tau(\varphi),V)v_g=0,
\end{equation}
where $V\in \Gamma(\varphi^{-1}TN)$ is a variation vector field along $\varphi$ which is given by 
$V(x)=\frac{d}{dt}\big\vert_{t=0}\varphi_t(x)\in T_{\varphi(x)}N$, 
$(x\in M)$, 
and  the {\em tension field} is given by 
$\tau(\varphi)
=\sum_{i=1}^mB(\varphi)(e_i,e_i)\in \Gamma(\varphi^{-1}TN)$, 
where 
$\{e_i\}_{i=1}^m$ is a locally defined orthonormal frame field on $(M,g)$, 
and $B(\varphi)$ is the second fundamental form of $\varphi$ 
defined by 
\begin{align}
B(\varphi)(X,Y)&=(\widetilde{\nabla}d\varphi)(X,Y)\nonumber\\
&=(\widetilde{\nabla}_Xd\varphi)(Y)\nonumber\\
&=\overline{\nabla}_X(d\varphi(Y))-d\varphi(\nabla_XY),
\end{align}
for all vector fields $X, Y\in {\frak X}(M)$. 
Here, 
$\nabla$, and
$\nabla^h$, 
 are Levi-Civita connections on $TM$, $TN$  of $(M,g)$, $(N,h)$, respectively, and 
$\overline{\nabla}$, and $\widetilde{\nabla}$ are the induced ones on $\varphi^{-1}TN$, and $T^{\ast}M\otimes \varphi^{-1}TN$, respectively. By (2.1), $\varphi$ is {\em harmonic} if and only if $\tau(\varphi)=0$. 
\par
The second variation formula is given as follows. Assume that 
$\varphi$ is harmonic. 
Then, 
\begin{equation}
\frac{d^2}{dt^2}\bigg\vert_{t=0}E(\varphi_t)
=\int_Mh(J(V),V)v_g, 
\end{equation}
where 
$J$ is an elliptic differential operator, called the
{\em Jacobi operator}  acting on 
$\Gamma(\varphi^{-1}TN)$ given by 
\begin{equation}
J(V)=\overline{\Delta}V-{\mathcal R}(V),
\end{equation}
where 
$\overline{\Delta}V=\overline{\nabla}^{\ast}\overline{\nabla}V
=-\sum_{i=1}^m\{
\overline{\nabla}_{e_i}\overline{\nabla}_{e_i}V-\overline{\nabla}_{\nabla_{e_i}e_i}V
\}$ 
is the {\em rough Laplacian} and 
${\mathcal R}$ is a linear operator on $\Gamma(\varphi^{-1}TN)$
given by 
${\mathcal R}(V)=
\sum_{i=1}^mR^N(V,d\varphi(e_i))d\varphi(e_i)$,
and $R^N$ is the curvature tensor of $(N,h)$ given by 
$R^h(U,V)=\nabla^h{}_U\nabla^h{}_V-\nabla^h{}_V\nabla^h{}_U-\nabla^h{}_{[U,V]}$ for $U,\,V\in {\frak X}(N)$.   
\par
J. Eells and L. Lemaire \cite{EL1} proposed polyharmonic ($k$-harmonic) maps and 
Jiang \cite{J} studied the first and second variation formulas of biharmonic maps. Let us consider the {\em bienergy functional} 
defined by 
\begin{equation}
E_2(\varphi)=\frac12\int_M\vert\tau(\varphi)\vert ^2v_g, 
\end{equation}
where 
$\vert V\vert^2=h(V,V)$, $V\in \Gamma(\varphi^{-1}TN)$.  
\par
The first variation formula of the bienergy functional 
is given by
\begin{equation}
\frac{d}{dt}\bigg\vert_{t=0}E_2(\varphi_t)
=-\int_Mh(\tau_2(\varphi),V)v_g.
\end{equation}
Here, 
\begin{equation}
\tau_2(\varphi)
:=J(\tau(\varphi))=\overline{\Delta}(\tau(\varphi))-{\mathcal R}(\tau(\varphi)),
\end{equation}
which is called the {\em bitension field} of $\varphi$, and 
$J$ is given in $(2.4)$.  
\par
A smooth map $\varphi$ of $(M,g)$ into $(N,h)$ is said to be 
{\em biharmonic} if 
$\tau_2(\varphi)=0$. 
By definition, every harmonic map is biharmonic. 
We say, for an immersion $\varphi:\,(M,g)\rightarrow (N,h)$  to be {\em proper biharmonic} if 
it is biharmonic but not harmonic (minimal). 
\vskip0.3cm\par
\subsection{Riemannian submersions.} 
We prepare with several notions on the Riemannian submersions. 
A $C^{\infty}$ mapping $\pi$ of a $C^{\infty}$ Riemannian manifold $(P,g)$ into another $C^{\infty}$ Riemannian manifold $(M,h)$ is called a {\em Riemannia submersion} if 
$(0)$ $\pi$ is surjective, $(1)$ the differential 
$d\pi=\pi_{\ast}:\,\,T_uP\rightarrow T_{\pi(u)}M\,\,(u\in P)$ 
of $\pi:\,\,P\rightarrow M$ 
is surjective for each $u\in P$, and 
$(2)$ each tangent space $T_uP$ at $u\in P$ has 
the direct decomposition: 
$$
T_uP={\mathcal V}_u\oplus {\mathcal H}_u,\qquad (u\in P), 
$$
which is orthogonal decomposition with respect to $g$ such  that ${\mathcal V}={\mbox{\rm Ker}}(\pi_{\ast\,u})\subset T_uP$ and 
$(3)$ the restriction of the differential $\pi_{\ast}=d\pi_u$ to 
${\mathcal H}_u$ is a surjective isometry,  
$\pi_{\ast}:\,\,({\mathcal H}_u,g_u)\rightarrow (T_{\pi(u)}M, h_{\pi(u)})$ for each $u\in P$ (cf. \cite{BW}, \cite{B}). 
A manifold $P$ is the total space of a Riemannian submersion over $M$ with the  projection $\pi:\,\,P\rightarrow M$ onto $M$, 
where $p=\dim P=k+m$, $m=\dim M$, and $k=\dim\pi^{-1}(x)$, $(x\in M)$. 
A Riemannian metric $g$ on $P$, called 
{\em adapted metric } on $P$ which satisfies  
\begin{equation}
g=\pi^{\ast}h+k
\end{equation}
where $k$ is the Riemannian metric on each fiber 
$\pi^{-1}(x)$, $(x\in M)$. Then, 
$T_uP$ has the orthogonal direct decomposition 
of the tangent space $T_uP$, 
\begin{equation}
T_uP={\mathcal V}_u\oplus {\mathcal H}_u,\qquad\quad u\in P, 
\end{equation}
where the subspace ${\mathcal V}_u=\mbox{\rm Ker}(\pi_{\ast}{}_u)$ at $u\in P$, 
the {\em vertical subspace}, and the subspace 
${\mathcal H}_u$ of $P_u$ is called 
{\em horizontal subspace} at $u\in P$ which is the orthogonal complement of ${\mathcal V}_u$ in $T_uP$ with respect to $g$. 
\par 
In the following, we fix a locally defined orthonormal frame field, called {\em adapted local orthonormal frame field} to the projection $\pi:\,\,P\rightarrow M$, 
$\{e_i\}_{i=1}^p$ corresponding to $(2.9)$ 
in such a way that 
\par $\bullet$
$\{e_i\}_{i=1}^m$ is a locally defined orthonormal basis of the horizontal 
\par\quad  subspace ${\mathcal H}_u$ $(u\in P)$, 
and 
\par $\bullet$ 
$\{ e_i\}_{i=1}^k$ 
is a locally defined orthonormal basis 
of the vertical \par\quad 
subspace ${\mathcal V}_u$ $(u\in P)$. 
\medskip\par
Corresponding to the decomposition $(2.9)$, the tangent vectors $X_u$, and $Y_u$ in $T_uP$ 
 can be defined by 
\begin{align}
&X_u=X_u^{{\rm V}}+X_u^{{\rm H}},\quad 
Y_u=Y_u^{{\rm V}}+Y_u^{{\rm H}},\\
&X_u^{{\rm V}},\,\,
Y_u^{{\rm V}}
\in {\mathcal V}_u,\quad 
X_u^{{\rm H}},\,\,
Y_u^{{\rm H}}\in {\mathcal H}_u
\end{align}
for $u\in P$. 
\par
Then, there exist a unique decomposition 
such that 
$$
g(X_u,Y_u)=h(\pi_{\ast}X_u,\pi_{\ast}Y_u)
+k(X^{\rm V}_u,Y^{\rm V}_u),
\quad X_u,\,\,Y_u\in T_uP,\,\,u\in P. 
$$
\vskip0.3cm\par
Then, let us recall 
the following definitions for our question:
 \vskip0.3cm\par
 {\bf Definition 2.1.}
 (1) The projection $\pi:\,(P,g)\rightarrow (M,h)$ is 
 to be {\em harmonic} if 
 the tension field vanishes, $\tau(\pi)=0$, and\vskip0.1cm\par
 (2) the projection $\pi:\,(P,g)\rightarrow (M,h)$ 
 is to be {\em biharmonic} if, the bitension field vanishes, 
 $\tau_2(\pi)=J(\tau(\pi))=0$.
\vskip0.3cm\par\noindent
 We define the Jacobi operator  $J$ for the projection $\pi$ by 
 \begin{equation}
 J(V):=\overline{\Delta}V-{\mathcal R}(V),\qquad V\in \Gamma(\pi^{-1}TM).
 \end{equation}
 Here,  
 \begin{align}
 \overline{\Delta}V&:= 
 -\sum_{i=1}^p
 \left\{\overline{\nabla}_{e_i}(
 {\overline{\nabla}}_{e_i}V)
 -{\overline{\nabla}}_{{\nabla}_{e_i}e_i}V\right\}
 =\overline{\Delta}_{\mathcal H}V+\overline{\Delta}_{\mathcal V}V.
 \end{align}
 where 
 \begin{align}
 \overline{\Delta}_{\mathcal H}V
 &=-\sum_{i=1}^m
 \left\{{\overline{\nabla}}_{e_i}
 ({\overline{\nabla}_{e_i}}V)
 -{\overline{\nabla}}_{{\nabla_{e_i}}e_i}V\right\},\\
 \overline{\Delta}_{\mathcal V}V
 &=-\sum_{i=1}^k
 \big\{{\overline{\nabla}}_{A^{\ast}_{m+i}}
 (\overline{\nabla}_{A^{\ast}_{m+i}}V) 
 -{\overline{\nabla}}_
{\nabla_{A^{\ast}_{m+i}}A^{\ast}_{m+i}}V\big\},
\end{align}
 \par \noindent
for $V\in \Gamma(\pi^{-1}TM)$, respectively. 
 Recall, 
 $\{e_i\}_{i=1}^p$ is a local orthonormal frame 
 field on $(P,g)$, 
 $\{e_i\}_{i=1}^m$ is a local orthonormal horizontal 
 field on $(M,h)$ and 
 $\{e_{m+i,\,u}\}_{i=1}^k$ $(u\in P)$ is an orthonormal 
 frame field on the vertical space 
 ${\mathcal V}_u$ 
 ($u\in P$). 
 We call 
 $\overline{\Delta}_{\mathcal H}$, the {\em horizontal Laplacian}, and $\overline{\Delta}_{\mathcal V}$, the {\em vertical Laplacian}, respectively. 
\section{The reduction of the biharmonic equation}
\subsection{}
Hereafter, we treat with the above problem more precisely in the case $\dim(\pi^{-1}(x))=1,\,\,(u\in P,\,\,\pi(u)=x)$. 
Let $\{e_1,\,\,e_1,\,\,\ldots,\,\,e_m\}$ be an adapted local orthonormal frame field being $e_{n+1}=e_m$, vertical. 
The frame fields 
$\{ e_i:\,\,i=1,2,\ldots,n\}$ are 
the basic orthonormal frame field 
on $(P,g)$ 
corresponds to an orthonormal 
frame field $\{\epsilon_1,\,\,\epsilon_2,\,\,\ldots,\,\,\epsilon_n\}$ on $(M,g)$.  
Here, a vector field $Z\in {\frak X}(P)$ is {\em basic} if 
$Z$ is horizontal and $\pi$-related to a vector field 
$X\in {\frak X}(M)$. 
\par
In this section, we determine the biharmonic equation 
precisely in the case that 
$p=m+1=\dim P$, $m=\dim M$, and $k=\dim \pi^{-1}(x)=1$ 
$(x\in M)$. 
Since $[V,Z]$ is a vertical field on $P$ if $Z$ is basic and $V$ is vertical (cf. \cite{ON}, p. 461).  Therefore, 
for each $i=1,\ldots,n$, 
$[e_i,e_{n+1}]$ is vertical, so we can write as follows. 
\begin{equation}
[e_i,e_{n+1}]=\kappa_i\,e_{n+1}, \quad
i=1,\,\,\ldots,\,\,n
\end{equation}
where $\kappa_i\in C^{\infty}(P)$ ($i=1,\ldots,n$). 
For two vector fields 
$X,\,\,Y$ on $M$, let $X^{\ast},\,\,Y^{\ast}$, be the 
horizontal vector fields on $P$. 
Then, 
$[X^{\ast},Y^{\ast}]$ is a vector field on $P$ which is 
$\pi$-related to a vector field $[X,Y]$ on $M$ (for instance, 
\cite{U2}, p. 143).  Thus, for $i,\,\,j=1,\ldots,n$, 
$[e_i,e_j]$ is $\pi$-related to 
$[\epsilon_i,\epsilon_j]$, and we may write as 
\begin{equation}
[e_i,e_j]=\sum_{k=1}^{n+1}D^k_{ij}\,e_k,
\end{equation}
where $D^k_{ij}\in C^{\infty}(P)\,\,
(1\leq i,\,j\leq n;\,1\leq k\leq n+1)$.
\subsection{The tension field} 
In this subsection, we calculate the tension 
field $\tau(\pi)$.  
We show that 
\begin{equation}
\tau(\pi)=-d\pi\left(\nabla_{e_{n+1}}e_{n+1}
\right)=-\sum_{i=1}^n\kappa_i\,\epsilon_i.
\end{equation}
Indeed, we have 
\begin{align}
\tau(\pi)&=\sum_{i=1}^m\left\{\nabla^{\pi}_{e_i}d\pi(e_i)-d\pi\left(\nabla_{e_i}e_i
\right)
\right\}\nonumber\\
&=\sum_{i=1}^n\left\{\nabla^{\pi}_{e_i}d\pi(e_i)-d\pi\left(\nabla_{e_i}e_i
\right)
\right\}
+\nabla^{\pi}_{e_{n+1}}d\pi(e_{n+1})-d\pi\left(\nabla_{e_{n+1}}e_{n+1}
\right)\nonumber\\
&=-d\pi\left(\nabla_{e_{n+1}}e_{n+1}
\right)\nonumber\\
&=-\sum_{i=1}^n\kappa_i\epsilon_i.\nonumber
\end{align}
Because, 
for $i,\,j=1,\ldots,n$, $d\pi(\nabla_{e_i}e_j)=\nabla^h_{\epsilon_i}\epsilon_j$, and 
$\nabla^{\pi}_{e_i}d\pi(e_i)=\nabla^h_{d\pi(e_i)}d\pi(e_i)=\nabla^h_{\epsilon}\epsilon_i$. Thus, we have 
\begin{equation}
\sum_{i=1}^n\left\{
\nabla^{\pi}_{e_i}d\pi(e_i)-d\pi\left(\nabla_{e_i}e_i
\right)
\right\}=0.
\end{equation}
Since $e_{n+1}=e_m$ is vertical, $d\pi(e_{n+1})=0$, so that 
$\nabla^{\pi}_{e_{n+1}}d\pi(e_{n+1})=0$. 
\par
Furthermore, we have, 
by definition of the Levi-Chivita connection, we have, 
for $i=1,\ldots,n$, 
$$
2g(\nabla_{e_{n+1}e_{n+1}},e_i)=2g(e_{n+1},[e_i,e_{n+1}])=2\kappa_{i},
$$
and $2g(\nabla_{e_{n+1}}e_{n+1},e_{n+1})=0$. Therefore, we have 
$$
\nabla_{e_{n+1}}e_{n+1}=\sum_{i=1}^n\kappa_ie_i,
$$
and then, 
\begin{equation}
d\pi\left(\nabla_{e_{n+1}}e_{n+1}
\right)=\sum_{i=1}^n\kappa_i\epsilon_i.
\end{equation}
Thus, we obtain (3.3). \qed
\subsection{The bitension field} 
Let us recall first the bitension field $\tau_2(\pi)$ is given by 
\begin{align}
\tau_2(\pi)&=-\sum_{i=1}^m\left\{
\nabla^{\pi}_{e_i}\left(\nabla^{\pi}_{e_i}\tau(\pi)
\right)
-\nabla^{\pi}_{\nabla_{e_i}e_i}\tau(\pi)
\right\}\nonumber\\
&\quad -\sum_{i=1}^m
R^h(\tau(\pi),d\pi(e_i))d\pi(e_i).
\end{align}
\par
First, since $d\pi(e_i)=\epsilon_{i},\,\, i=1,\ldots,n$, 
we have 
\begin{align}
\sum_{i=1}^nR^h(\tau(\pi),d\pi(e_i))d\pi(e_i)&=\sum_{i=1}^nR^h(\tau(\pi),\epsilon_i)\epsilon_i\nonumber\\
&={\mbox {\rm Ric}}^h(\tau(\pi)).
\end{align}
On the other hand, we calculate the first term of (3.6) 
for $\tau_2(\pi)$. 
\medskip\par
(The first step)\,\,To calculate $\nabla^{\pi}_{e_i}\tau(\pi)$ $(i=1,\,\ldots,\,m=n+1)$, we want to show
\begin{equation}
\nabla^{\pi}_{e_i}\tau(\pi)=
\left\{\begin{aligned}
&-\sum_{j=1}^n\left\{(e_i\kappa_j)\epsilon_j+\kappa_j\,\nabla^h_{\epsilon_i}\epsilon_j
\right\}\quad (i=1,\ldots,n),\\
&\qquad 0\qquad\qquad\qquad\qquad\qquad\,\,
(i=n+1).
\end{aligned}
\right.
\end{equation}   
\par
Because, if $i=1,\ldots,n$, by noticing $\kappa_j\in C^{\infty}(P)$, $(j=1,\ldots,n))$, we have by $(3.3)$, 
\begin{align}
\nabla^{\pi}_{e_i}\tau(\pi)&=\nabla^{\pi}_{e_i}\left(
-\sum_{j=1}^n\kappa_j\,\epsilon_j
\right)\nonumber\\
&=-\sum_{j=1}^n\left\{(e_i\,\kappa_j)\,\epsilon_j+\kappa_j\,\nabla^{\pi}_{e_i}\epsilon_j
\right\}\nonumber\\
&=-\sum_{j=1}^n\left\{(e_i\,\kappa_j)\,\epsilon_j+\kappa_j\,\nabla^h_{\epsilon_i}\epsilon_j
\right\}, 
\end{align}
since $\nabla^{\pi}_{e_i}\epsilon_j=\nabla^h_{d\pi(e_i)}\epsilon_j=\nabla^h_{\epsilon_i}\epsilon_j$. 
Furthermore, for $i=n+1$, we have 
\begin{equation}
\nabla^{\pi}_{e_{n+1}}\tau(\pi)=\nabla^h_{d\pi(e_{n+1})}\tau(\pi)=0. 
\end{equation}
To show (3.10), recalling the definition of the parallel displacement of the connection, 
let
$P_{\pi\circ\sigma(t)}:\,\,T_{\pi(\sigma(0))}M\rightarrow T_{\pi(\sigma(t))}M$ be the parallel transport 
with respect to $(M,h)$ along a smooth curve  in $P$.  
Then, since 
$\sigma(t)\in P,\,\,\epsilon<t<\epsilon$ with $\sigma(0)=x\in P$ and 
$\dot{\sigma}(0)=e_{n+1\,\,x}\in T_xP$, 
for every $V\in \Gamma(\pi^{-1}TM)$, and then, 
\begin{align}
\nabla^{\pi}_{e_{n+1}}V(x)=\frac{d}{dt}\bigg\vert_{t=0}P^{-1}_{\pi\circ\sigma(t)}V(\sigma(t))
=\frac{d}{dt}\bigg\vert_{t=0}P^{-1}_{\pi(x)}V(\sigma(t))=0,
\end{align} 
since 
$\pi(\sigma(t))=\pi(\sigma(0))=\pi(x)\in P$ because $e_{n+1}$ is a vertical vector field of the Riemannian submersion $\pi:\,\,(P,g)\rightarrow (M,h)$.  
\medskip\par
(The second step) To calculate $\nabla^{\pi}_{\nabla_{e_i}e_i}\tau(\pi)$ $(i=1,\,\ldots,\,m=n+1)$, we have 
\begin{equation}
\nabla^{\pi}_{\nabla_{e_i}e_i}\tau(\pi)=
\left\{\begin{aligned}
&-\sum_{j=1}^n\left\{(\nabla_{e_i}e_i\,\kappa_j)\epsilon_j+\kappa_j\,\nabla^h_{\nabla^h_{\epsilon_i}\epsilon_i}\epsilon_j
\right\} (i=1,\ldots,n),\\
&-\sum_{\ell,\,j=1}^n\left\{\kappa_{\ell}\,(e_{\ell}\kappa_j)\,\epsilon_j+\kappa_{\ell}\kappa_j\,\nabla^h_{\epsilon_{\ell}}\epsilon_j\right\}
 (i=n+1).
\end{aligned}
\right.
\end{equation}   
\par
Indeed, for a vector field $\nabla_{e_i}e_i$ on $P$ ($i=1,\ldots,n$), we only have to see that 
\begin{equation}
d\pi(\nabla_{e_i}e_i)=\nabla^h_{\epsilon_i}\epsilon_i, 
\end{equation} 
which yields the first equation of $(3.12)$. 
To see (3.13), 
we have to see the following equations:  
\begin{align}
\nabla_{e_i}e_i&={\mathcal V}(\nabla_{e_i}e_i)+{\mathcal H}(\nabla_{e_i}e_i)\nonumber\\
&=A_{e_i}e_i+{\mathcal H}(\nabla_{e_i}e_i)\quad (\mbox{cf. the fourth of Lemma 3 in \cite{ON}, p. 461})\nonumber\\
&=\frac12{\mathcal V}[e_i,e_i]+ {\mathcal H}(\nabla_{e_i}e_i) \qquad (\mbox{cf. Lemma 2 in \cite{ON}, p. 461)}\nonumber\\
&={\mathcal H}(\nabla_{e_i}e_i). 
\end{align}
Here, since ${\mathcal H}(\nabla_{e_i}e_i)$ is a basic vector field 
corresponding to 
$\nabla^h_{\epsilon_i}\epsilon_i$ (cf. the third of Lemma 1 in \cite{ON}, p. 460), we have 
$d\pi(\nabla_{e_i}e_i)=d\pi({\mathcal H}(\nabla_{e_i}e_i))=\nabla^h_{\epsilon_i}\epsilon_i$, i.e.,  $(3.13)$. 
Then, we have 
\begin{align}
\nabla^{\pi}_{\nabla_{e_i}e_i}\tau(\pi)&=\sum_{\nabla_{e_i}e_i}(-\sum_{j=1}^n\kappa_j\epsilon_j)\nonumber\\
&=-\sum_{j=1}^n\left\{
(\nabla_{e_i}e_i\,\kappa_j)\,\epsilon_j+\kappa_j\nabla^{\pi}_{\nabla_{e_i}e_i}\epsilon_j
\right\}\nonumber\\
&=-\sum_{j=1}^n\left\{
(\nabla_{e_i}e_i\,\kappa_j)\,\epsilon_j+\kappa_j\,\nabla^h_{\nabla^h_{\epsilon_i}\epsilon_i}\epsilon_j
\right\}, 
\end{align}
which is the first equation of (3.12). 
To see the second equation of $(3.12)$, 
recall (3.5)
$d\pi(\nabla_{e_{n+1}}e_{n+1})=\sum_{i=1}^n\kappa_i\epsilon_i$  and also the first equation of $(3.8)$. Then, 
we have 
\begin{align}
\nabla^{\pi}_{\nabla_{e_{n+1}e_{n+1}}}\tau(\pi)&=
-\nabla^h_{\left(\sum_{i=1}^n\kappa_i\epsilon_i\right)}\sum_{j=1}^n\kappa_j\epsilon_j\nonumber\\
&=-\sum_{i,j=1}^n\left\{
\kappa_i\,\epsilon_i(\kappa_j)\,\epsilon_j+\kappa_i\kappa_j\,\nabla^h_{\epsilon_i}\epsilon_j
\right\}, 
\end{align}
which implies the second equation of $(3.12)$. 
\medskip\par
(The third step) We calculate $\nabla^{\pi}_{e_i}(\nabla^{\pi}_{e_i}\tau(\pi))$. Indeed, we have 
\begin{align}
&\nabla^{\pi}_{e_i}\left(\nabla^{\pi}_{e_i}\tau(\pi)\right)=\nabla^{\pi}_{e_i}\left(
-\sum_{j=1}^n\left\{(e_i\kappa_j)\,\epsilon_j+\kappa_j\,\nabla^h_{\epsilon_i}\epsilon_j
\right\}
\right)\nonumber\\
&=-\sum_{j=1}^n\left\{
e_i(e_i\kappa_j)\,\epsilon_j+(e_i\kappa_j)\nabla^{\pi}_{e_i}\epsilon_j+(e_i\kappa_j)\,\nabla^h_{\epsilon_i}\epsilon_j
+\kappa_j\nabla^{\pi}_{e_i}\left(\nabla^h_{\epsilon_i}\epsilon_j
\right)
\right\}, 
\end{align}
where 
\begin{equation}
\left\{\begin{aligned}
&\nabla^{\pi}_{e_i}\epsilon_j=\nabla^h_{d\pi(e_i)}\epsilon_j=\nabla^h_{\epsilon_i}\epsilon_j,\\
&\nabla^{\pi}_{e_i}(\nabla^h_{\epsilon_i}\epsilon_j)=\nabla^h_{d\pi(e_i)}(\nabla^h_{\epsilon_i}\epsilon_j)=\nabla^h_{\epsilon_i}(\nabla^h_{\epsilon_i}\epsilon_j).
\end{aligned}\right.
\end{equation}
Then we have, for $i=1,\ldots,n$, 
\begin{equation}
\left\{
\begin{aligned}
&\nabla^{\pi}_{e_i}\left(\nabla^{\pi}_{e_i}\tau(\pi)\right)=
-\sum_{j=1}^n\left\{
e_i(e_i\kappa_j)\,\epsilon_j+2(e_i\kappa_j)\,\nabla^h_{\epsilon_i}\epsilon_j+\kappa_j\,\nabla^h_{\epsilon_i}(\nabla^h_{\epsilon_i}\epsilon_j)
\right\}, \\
&\nabla^{\pi}_{e_{n+1}}
(\nabla^{\pi}_{e_{n+1}}\tau(\pi))=0,
\\
&\nabla^{\pi}_{\nabla_{e_i}e_i}\tau(\pi)
=-\sum_{j=1}^n
\left\{
(\nabla_{e_i}e_i\,\kappa_j)\,\epsilon_j+\kappa_j\,\nabla^h_{\nabla^h_{\epsilon_i}\epsilon_i}\epsilon_j
\right\},\\
&\nabla^{\pi}_{\nabla_{e_{n+1}}e_{n+1}}\tau(\pi)
=-\sum_{i,j=1}^n\left\{\kappa_i(e_i\kappa_j)\,\epsilon_j+\kappa_i\kappa_j\,\nabla^h_{\epsilon_i}\epsilon_j\right\}.
\end{aligned}
\right.
\end{equation}
\medskip\par
In the following, by using these formulas, we can proceed to calculate 
\begin{align*}
\tau_2(\pi)&=\overline{\Delta}^h\,\tau(\pi)-{\mbox{Ric}}^h(\tau(\pi))\nonumber\\
&=-\sum_{i=1}^m\left\{\nabla^{\pi}_{e_i}\left(\nabla^{\pi}_{e_i}\tau(\pi)
\right)-\nabla^{\pi}_{\nabla_{e_i}e_i}\tau(\pi)
\right\}-{\mbox{Ric}}^h(\tau(\pi)). 
\end{align*}
(The fourth step)  Indeed, we have 
\begin{align}
\tau_2(\pi)&=\overline{\Delta}^h\,\tau(\pi)-{\mbox{Ric}}^h(\tau(\pi))\nonumber\\
&=-\sum_{i=1}^m\left\{\nabla^{\pi}_{e_i}\left(\nabla^{\pi}_{e_i}\tau(\pi)
\right)-\nabla^{\pi}_{\nabla_{e_i}e_i}\tau(\pi)
\right\}-{\mbox{Ric}}^h(\tau(\pi))\nonumber\\
&=\sum_{i,j=1}^{n}\bigg\{
e_i(e_i\kappa_j)\,\epsilon_j+2(e_i\kappa_j)\nabla^h_{\epsilon_i}\epsilon_j+\kappa_j\,\nabla^h_{\epsilon_i}(\nabla^h_{\epsilon_i}\epsilon_j)\nonumber\\
&\qquad \qquad 
-(\nabla_{e_i}e_i\,\kappa_j)\epsilon_j-\kappa_j\,\nabla^h_{\nabla^h_{\epsilon_i}\epsilon_i}\epsilon_j
-\kappa_i(e_i\kappa_j)\,\epsilon_j-\kappa_i\kappa_j\,\nabla^h_{\epsilon_i}\epsilon_j
\bigg\}\nonumber\\
&\qquad+{\mbox{Ric}}^h\bigg(\sum_{j=1}^n\kappa_j\epsilon_j\bigg)\nonumber
\\
&=\sum_{j=1}^n\sum_{i=1}^n\bigg\{e_i(e_i\,\kappa_j)-\nabla_{e_i}e_i\,\kappa_j
\bigg\}\epsilon_j
+2\sum_{j=1}^n\nabla^h_{(\sum_{i=1}^n(e_i\kappa_j)\epsilon_i)}\epsilon_j\nonumber\\
&\qquad+\sum_{j=1}^n\kappa_j\sum_{i=1}^n\bigg\{
\nabla^h_{\epsilon_i}\nabla^h_{\epsilon_i}\epsilon_j-\nabla^h_{\nabla^h_{\epsilon_i\epsilon_i}}\epsilon_j
\bigg\}
-\nabla_{(\sum_{i=1}^n\kappa_i\epsilon_i)}\sum_{j=1}^n\kappa_j\epsilon_j\nonumber\\
&\qquad +{\mbox{Ric}}^h(\sum_{j=1}^n\kappa_j\epsilon_j)
\nonumber
\end{align}
\begin{align}
&=\sum_{j=1}^n\bigg\{-\Delta\kappa_j-e_{n+1}(e_{n+1}\kappa_j)+\nabla_{e_{n+1}}e_{n+1}\,\kappa_j\bigg\}\,\epsilon_j\nonumber\\\
&\qquad+2\sum_{j=1}^n\nabla^h_{(\sum_{i=1}^n(e_i\kappa_j)\epsilon_i)}\epsilon_j
-\sum_{j=1}^n\kappa_j(\overline{\Delta}^h\epsilon_j) 
-\nabla_{(\sum_{i=1}^n\kappa_i\epsilon_i)}\sum_{j=1}^n\kappa_j\epsilon_j \nonumber\\
&\qquad+{\mbox{Ric}}^h(\sum_{j=1}^n\kappa_j\epsilon_j)\nonumber\\
&=\sum_{j=1}^n(-\Delta^h\kappa_j)\,\epsilon_j\nonumber\\
&\qquad +2\sum_{j=1}^n\nabla^h_{(\sum_{i=1}^n(e_i\kappa_j)\epsilon_i)}\epsilon_j
-\sum_{j=1}^n\kappa_j(\overline{\Delta}^h\epsilon_j) 
-\nabla_{(\sum_{i=1}^n\kappa_i\epsilon_i)}\sum_{j=1}^n\kappa_j\epsilon_j \nonumber\\
&\qquad+{\mbox{Ric}}^h(\sum_{j=1}^n\kappa_j\epsilon_j).
\end{align}
Since 
\begin{align}
\overline{\Delta}^h(\kappa_j\epsilon_j)=
(\overline{\Delta}^h\kappa_j)\epsilon_j-2\sum_{i=1}^n(e_i\kappa_j)\,\nabla^h_{\epsilon_i}\epsilon_j+\kappa_j(\overline{\Delta}^h\epsilon_j), 
\end{align}
we obtain 
\begin{align}
\tau_2(\pi)=-\Delta^h\big(\sum_{j=1}^n\kappa_j\epsilon_j\big)
-\nabla^h_{(\sum_{i=1}^n\kappa_i\epsilon_i)}\sum_{j=1}^n\kappa_j\epsilon_j+{\mbox{Ric}}^h\big(\sum_{j=1}^n\kappa_j\epsilon_j\big).
\end{align}
Thus, we obtain the following theorem: 
\medskip\par
\begin{th} Let $\pi:\,\,(P,g)\rightarrow (M,h)$ be a Riemannian submersion over $(M,h)$. Then, 
\par
$(1)$ The tension field $\tau(\pi)$ of $\pi$ is given by 
\begin{align}
\tau(\pi)=-\sum_{i=1}^n\kappa_i\epsilon_i,
\end{align}
where $\kappa_i\in C^{\infty}(P)$, $(i=1,\ldots,n)$. 
\par
$(2)$ The bitension field $\tau_2(\pi)$ of $\pi$ is given by 
\begin{align}
\tau_2(\pi)=-\overline{\Delta}^h\big(\sum_{j=1}^n\kappa_j\epsilon_j\big)+\nabla^h_{(\sum_{i=1}^n\kappa_i\epsilon_i)}\sum_{j=1}^n\kappa_j\epsilon_j+
{\mbox{\rm Ric}}^h\big(\sum_{j=1}^n\kappa_j\epsilon_j\big).
\end{align}
\end{th}
\medskip\par
\begin{rem}  
The bitension field $\tau_2(\pi)$ for $\pi$ has been 
obtained in a different way 
by Akyol and Ou \cite{AO} 
in which has referenced our paper. 
\end{rem}
\medskip\par
\begin{prop}
Let $\pi:\,\,(P,g)\rightarrow (M,h)$ be a Riemannian submersion whose base manifold 
$(M,h)$ has non-positive Ricci curvature. 
Assume that $\pi:,\,(P,g)\rightarrow (M,h)$ is biharmonic. 
 Then the tension field $X:=\tau(\pi)$ is parallel, i.e., 
$\nabla^hX=0$ if we assume $\mbox{div}(X)=0$. 
\end{prop}
{\em Proof} \quad 
Assume that $\pi:\,\,(P,g)\rightarrow (M,h)$ is biharmonic, 
i.e., 
$$
0=\tau_2(\pi)=-{\overline{\Delta}}^hX-\nabla^h_XX+\mbox{\rm Ric}^h(X).
$$
Then, we have 
\begin{align}
0&\leq \int_M{\overline{\nabla}}^hX, 
{\overline{X}}^hX)\,v_h\nonumber\\
&=\int_Mh({\overline{\Delta}}^hX,X)\,g_h\nonumber\\
&=-\int_Mh(\nabla^h_XX,X)\,v_h
+\int_Mh({\mbox{\rm Ric}}^h(X),X)\,v_h\nonumber\\
&=-\frac12\int_MX\,\cdot\,h(X,X)\,v_h
+\int_Mh({\mbox{\rm Ric}}^h(X),X)\,v_h\nonumber\\
&=\int_M({\mbox{\rm Ric}}^h(X),X)\,v_h\leq 0.
\end{align}
The second equality from below holds, 
due to Gaffney's theorem (cf. Theorem 2.2 in \cite{OSU2}), 
$\int _MXf\,v_h=0$ $(f\in C^1(M))$ 
if $\mbox{\rm div}(X)=0$.   
The last inequality holds for non-positive Ricci curvature 
of $(M,h)$. 
Therefore, we have 
$$
0=h(\mbox{\rm Ric}(X),X)\,v_h=\int_Mh({\overline{X}}^h,
{\overline{X}}^h)\,v_h.
$$
Thus, we have ${\overline{\nabla}}^hX=0$. 
\qed
\vskip0.6cm\par
\section{Einstein manifolds}
\subsection{}
Regarding the orthogonal direct decomposition: 
\begin{equation}
{\frak X}(M)=\{X\in {\frak X}(M)\vert\,\,\mbox{\rm div}(X)=0\}
\oplus\{\nabla\,f\in {\frak X}(M)\vert\,\,f\in C^{\infty}(M)\},
\end{equation}
we obtain the following theorems: 
\begin{th}
Let $\pi:\,\,(P,g)\rightarrow (M,h)$ be 
a compact Riemannian submersion 
over a weakly stable Einstein manifold $(M,g)$ whose Ricci tensor 
$\rho^h$ satisfies 
$\rho^h=c\,\mbox{\rm Id}$ for some constant $c$. 
Assume that $\pi$ is biharmonic, i.e., 
\begin{equation}
\tau_2(\pi)=-\overline{\Delta}^hX+\nabla^h_XX+{\mbox{\rm Ric}}^h(X)=0, 
\end{equation}
where $X=\sum_{i=1}^n\kappa_i\epsilon_i$. 
Assume that $\mbox{\rm div}X=0$. 
Then, 
\begin{equation}
\left\{
\begin{aligned}
\overline{\Delta}^hX&=cX,\\
\nabla^h_XX&=0.
\end{aligned}
\right.
\end{equation}
\end{th}
{\em Proof}\quad 
Let $X=\sum_{i=1}^{\infty}X_i$ where 
$\Delta^HX_i=\lambda_iX_i$ satisfying that 
\par\noindent
$\int_Mh(X_i,X_j)\,v_h=\delta_{ij}$. $\Delta^H$ corresponds to the Laplacian $\Delta^1$ acting 
on the space 
$A^1(M)$ of $1$-forms on $(M,h)$. 
By $(4.1)$, 
\begin{align}
-\nabla^h_XX&=\overline{\Delta}^hX-cX\nonumber\\
&=\sum_{i=1}^{\infty}\lambda_iX_i-2c\sum_{i=1}^{\infty}X_i
\\
&=\sum_{i=1}^{\infty}(\lambda_i-2c)X_i
\end{align}
since $\Delta^H=\overline{\Delta}^h+\rho^h=\overline{\Delta}+c\,{\mbox{\rm Id}}$.  
Since ${\mbox{\rm div}}(X)=0$, 
\begin{align}
0&=-\frac12\int_MX\,\cdot\,h(X,X)\,v_h\nonumber\\
&=-\int_Mh(\nabla^h_XX,X)\,v_h\nonumber\\
&=\int_Mh(\sum_{i=1}^{\infty}(\lambda_i-2c)X_i,\sum_{j=1}^{\infty}X_j)\,v_h\nonumber\\
&=\sum_{i=1}^{\infty}(\lambda_i-2c).
\end{align}
If $(M,h)$ is weakly stable, i.e., 
$2c\leq \lambda^1_1(h)\leq \lambda_i\,\,(i=1,2,\ldots)$, 
then we have 
$$
\lambda_i=2c\qquad (i=1,2,\ldots). 
$$
Therefore, we have 
$$\overline{\Delta}^hX+cX=\Delta^HX=\sum_{i=1}^{\infty}\lambda_iX_i=2c\sum_{i=1}^{\infty}X_i=2cX. 
$$
Therefore, 
$$
\left\{\begin{aligned}
&\overline{\Delta}^hX=cX,\\
&\nabla^h_XX=0. 
\end{aligned}
\right.
$$
We have Theorem 4.1. 
\qed
\medskip\par
We have immediately the following theorem and corrollary: 
\begin{th} 
Let $\pi:\,\,(P,g)\rightarrow (M,h)$ be 
a compact Riemannian submersion 
over an irreducible compact Hermitian 
symmetric space $(M,h)=(K/H, h)$ where 
$K$ is a compact semi-simple Lie group, 
and $H$, a closed subgroup of $K$, 
$h$, an invariant Riemannan metric on $M=K/H$, 
respectively.   Let $X\in {\frak k}$ be an invariant vector field 
on $M$. 
Then, ${\mbox{\rm div}}X=0$, and that 
\begin{equation}
\left\{
\begin{aligned}
\overline{\Delta}^hX&=cX,\\
\nabla^h_XX&=0.
\end{aligned}
\right.
\end{equation}
\end{th}
\begin{cor} 
Let $\pi:\,\,(P,g)\rightarrow (M,h)$ be 
a principal $S^1$- bundle over 
an $n$-dimensional compact Hermitian symmetric space $(M,h)$. 
Then, 
\begin{equation}
\tau(\pi)=-\sum_{j=1}^n\kappa_j\widetilde{\epsilon}_j\in \Gamma(\pi^{-1}TM).
\end{equation}
If $X=\sum_{i=1}^n\kappa_j\epsilon_j$ is a non-vanishing Killing vector field on $(M,h)$, $\pi:\,\,(P,g)\rightarrow (M,h)$ is biharmonic, but not harmonic. 
\end{cor}
\medskip\par
\subsection{}
Regarding (4.1), we now consider the case 
$\{\nabla\,f\in {\frak X}(M)\vert\,\,f\in C^{\infty}(M)\}$. 
 Recall a theorem of  M. Obata 
on a compact K\"{a}hler-Einstein Riemannian manifold 
$(M,h)$
 (\cite{U2}, p. 181), 
the first non-zero positive eigenvalue 
$\lambda_1(h)$ of $(M,h)$  satisfies that 
\begin{equation}
\lambda_1(h)\geq 2c,
\end{equation}
and if the equality $\lambda_1(h)=2c$ holds, 
the corresponding eigenfunction $f$ with the eigenvalue $2c$ satisfies that 
$\nabla\, f$ is an analytic vector field on $M$ (\cite{U2}, p. 174) and 
\begin{equation}
J_{\mbox{\rm id}}(\nabla\,f)=0, 
\end{equation}
where $J_{\mbox{\rm id}}$ is the Jacobi operator given by 
$
J_{\mbox{\rm id}}:={\overline{\Delta}}^h-2\,\mbox{\rm Ric}.
$
\par
We apply the above to the our situation that 
$\pi:\,\,(P,g)\rightarrow (M,h)$ is a compact Riemannian submersion 
over a compact K\"{a}hler-Einstein manifold $(M,h)$ with 
$\mbox{\rm Ric}^h=c\,\mbox{Id}$, 
and assume that $\pi:\,\,(P,g)\rightarrow (M,h)$ is biharmonic, 
i.e., 
\begin{equation}
{\overline{\Delta}}^hX+{\nabla^h}_XX-\mbox{\rm Ric}^h(X)=0, 
\end{equation}
where $X=\tau(\pi)\in \Gamma(\pi^{-1}TM)$. 
\par
Thus, we can summarize the above as follows: 
\begin{th}
Assume that our $X=\tau(\pi)$ is of the form, 
$X=\nabla f$, where $f$ is the eigenfunction of the Laplacian 
$\Delta_h$ acting on $C^{\infty}(M)$ 
with the first eigenvalue $\lambda_1(h)=2c$. 
\par 
Then $X$ is an analytic vector field on $M$ (\cite{U2}, p. 174) and 
\begin{equation}
J_{\mbox{\rm id}}(X)=0, 
\end{equation}
where $J_{\mbox{\rm id}}$ is the Jacobi operator given by 
$
J_{\mbox{\rm id}}:={\overline{\Delta}}^h-2\,\mbox{\rm Ric}.
$
\par
Furthermore, we have 
\begin{equation}
\Delta_HX=2c\,X, \,\,\mbox{\rm i.e.}, \quad
{\overline{\Delta}}^hX=c\,X, 
\end{equation}
and also 
\begin{equation}
{\nabla^h}_XX=0.
\end{equation}
Here, $\Delta_H$ is the operator acting on 
   ${\frak X}(M)$ corresponding to the standard Laplacian 
   $\Delta:=d\delta+\delta d$ on the space $A^1(M)$ 
   of $1$-forms on $(M,h)$. 
   \end{th}
\medskip\par
\subsection{} 
In this part, we show 
\begin{prop}
Under the above situation, we have, at each point 
$p\in P$, 
\begin{equation}
\mbox{\rm div}(X)(p)
=\sum_{i=1}^ne_i\,\kappa_i\,(p),
\end{equation}
where 
$\{e_i\}_{i=1}^n$ is an orthonormal frame field on a neighborhood of each point $p\in P$ satisfying that
 $(\nabla_Y e_i)(p)=0$, $\forall\,\,Y\in T_pP$ $(i=1,\ldots,n)$.  
\end{prop}
{\em Proof.}\quad Let us recall 
$X:=\tau(\pi)=-\sum_{i=1}^n\kappa_i\,\widetilde{\epsilon}_i\in \Gamma(\pi^{-1}TM)$, 
where $\kappa_i\in C^{\infty}(P)$, $\widetilde{\epsilon}_i=\pi^{-1}\epsilon_i\in \Gamma(\pi^{-1}TM)$ defined by 
$$
\widetilde{\epsilon}_i(p):=(\pi^{-1}\epsilon_i)(p)
=\epsilon_{i\,\pi(p)},\quad (p\in P), 
$$ 
and $\{\epsilon_i\}_{i=1}^n$ is a locally defined orthonormal frame field on $(M,h)$. 
Here, note that, for $p\in P,\,\,\pi(p)=x\in M$,
$$
X(p)=-\sum_{i=1}^n\kappa_i(p)\widetilde{\epsilon}_i(p)=-\sum_{i=1}^n\kappa_i(p)\epsilon_i(x)\in T_xM. 
$$
Let $\widetilde{\nabla}$ be the induced connection 
on $\Gamma(\pi^{-1}TM)$ from the Levi-Civita connection 
$\nabla^h$ of $(M,h)$, and define ${\mbox{\rm div}}(X)\in C^{\infty}(P)$ by 
\begin{align}
{\mbox{\rm div}}(X)(p):&=\sum_{i=1}^m g_p(e_{i\,\,p},\,(\widetilde{\nabla}_{e_i}X)(p))
=\sum_{i=1}^mg_p(e_{i\,\,p},\nabla^h_{\pi_{\ast}e_i}X)\nonumber\\
&=\sum_{i=1}^ng_p(e_{i\,\,p},
(\widetilde{\nabla}_{e_i}X)(p)), 
\end{align}  
where $m=n+1=\dim (P)$. Because 
$\widetilde{\nabla}_{e_{n+1}}X(p)=0$ since, 
for a $C^1$ curve $\sigma$ in $P$ with $\sigma(0)=p,\,\,\sigma'(0)=(e_{n+1})_p\in T_pP$, we have $\pi\,\circ\,\sigma_t(s)=x,\,\,\forall \,\,0\leq s\leq t$. Therefore, 
we have 
$$
(\widetilde{\nabla}_{e_{n+1}}X)(x)=\nabla^h_{\pi_{\ast}e_{n+1}}X=\frac{d}{dt}\bigg\vert_{t=0}P^h_{\pi\,\circ\,\sigma_t}{}^{-1}X(\sigma(t))=0,
$$
where $P^h_{\pi\,\circ\,\sigma_t}:\,\,T_{\pi(p)}M\rightarrow T_{\pi(\sigma(t))}M$ is the parallel displacement along a $C^1$ curve $\pi\,\circ\,\sigma_t$ with respect to $\nabla^h$ on $(M,h)$. 
Then, for the RHS of (4.16), we have 
\begin{align}
{\mbox{\rm div}}(X)(p)&=\sum_{i=1}^ng_p(e_{i\,\,p},(\widetilde{\nabla}_{e_i}X)(p))\nonumber\\
&=\sum_{i=1}^ng_p(e_{i\,\,p},\widetilde{\nabla}_{e_i}(\sum_{j=1}^n\kappa_j\widetilde{\epsilon}_j))\nonumber\\
&=\sum_{i=1}^ng_p\big(e_{i\,\,p},\sum_{j=1}^n\big\{
e_i\kappa_j(p)\,\widetilde{\epsilon}_j(p)+\kappa_j(p)(\widetilde{\nabla}_{e_i}\widetilde{\epsilon}_j)(p)
\big\}\big)\nonumber\\
&=\sum_{i,j=1}^n(e_i\kappa_j)(p)\,g_p(e_{i\,\,p},\widetilde{\epsilon}_j(p))+\sum_{i,j=1}^n\kappa_j(p)\,g_p(e_{i\,\,p},(\widetilde{\nabla}_{e_i}\widetilde{\epsilon}_j)(p))\nonumber\\
&=
\sum_{i=1}^n(e_i\kappa_i)(p)-g_p(\sum_{i=1}^n\nabla^g_{e_i}e_i,\sum_{j=1}^n\kappa_j(p)\widetilde{\epsilon}_j)\nonumber\\
&=\sum_{i=1}^ne_i\kappa_i+g(\sum_{i=1}^n\nabla^g_{e_i}e_i,X),
\end{align}
since 
$$
g_p(e_{i\,\,p},(\widetilde{\nabla}_{e_i}\widetilde{\epsilon}_j)(p)=
e_{i\,\,p}g(e_i,\widetilde{\epsilon}_j)-g_p(\nabla^g_{e_i}e_i,\widetilde{\epsilon}_j)
=-g_p(\nabla^g_{e_i}e_i,\widetilde{\epsilon}_j)
$$
by means of $e_{i\,p}\,g(e_i,\widetilde{\epsilon}_j)=0$. 
By noticing that 
$g(\sum_{i=1}^n\nabla^g_{e_i}e_i,X)=0$ at the point $p\in P$ because of a choice of $\{e_i\}$, we obtain (4.15). 
\qed
\vskip0.6cm\par
\section{K\"{a}hler-Einstein flag manifolds} 
Let $(M,h)=(K/T,h)$ be a K\"{a}hler-Einstein flag manifold with 
$\mathrm{Ric}^h=c\,\,\mathrm{Id}$ for some $c>0$,  
where $T$ be a maximal torus in $K$, and let
$E_{\lambda}$, the line bundle over $K/T$ associated to 
 non-trivial homomorphism 
 $\lambda:\,\,T\rightarrow {\mathbb C}^{\ast}$. 
 Then, $E_{\lambda}$ is the totality of 
 all equivalence classes $[k,v]$ including 
 $(k,v)$ 
 with $k\in K$ and $v\in {\mathbb C}^{\ast}$ under the equivalence relation 
 $(k',v')\sim (k,v)$, i.e., 
 $k'=ka,\,\,v'=\lambda(a^{-1})v$ for some 
 $a\in T$. 
 Let 
 ${\mathcal S}_{\lambda}:=\{
 [k,u]\vert\,k\in K,\,\,u\in S^1
 \}=\{(k,u)\vert\,k\in K,\,\,u\in S^1\}\slash\sim$. Then, 
 ${\mathcal S}_{\lambda}$ is the circle bundle over 
 a flag manifold $K/T$ associted to $\lambda:\,\,T\rightarrow S^1$, where $S^1=\{u\in {\mathbb C}\vert\,\,\vert u\vert=1\}$. 
  Note that $m:=\dim {\mathcal S}=n+1$, 
  with $n=\dim M=\dim K/T$. 
\vskip0.3cm\par
{\bf Example 1.} \quad For $r=1,2,\ldots$, 
let 
\begin{align*}
K=SU(r+1)\supset T=
\bigg\{
&\begin{bmatrix}
e^{2\pi\sqrt{-1}\,\theta_1}&& {\mathrm O}\\
&\ddots&&\\
{\mathrm O}&&
e^{2\pi\sqrt{-1}\,\theta_{r+1}}\!\!\!\!\!\!
\end{bmatrix}
\big\vert\\
&
\,\,\theta_1,\,\,\ldots,\,\,\theta_{r+1}\in {\mathbb R},\,\,\theta_1+\cdots+\theta_{r+1}=0
\bigg\},
\end{align*}
and for ${\mathcal  I}=(a_1,\,\,\ldots,a_{r+1})\in {\mathbb Z}^{r+1}$, 
let 
$$
\lambda_{\mathcal I}:\,\,T\ni 
\begin{bmatrix}
e^{2\pi\sqrt{-1}\,\theta_1}&& {\mathrm O}\\
&\ddots&&\\
{\mathrm O}&&
e^{2\pi\sqrt{-1}\,\theta_{r+1}}\!\!\!\!\!\!
\end{bmatrix}
\mapsto 
e^{2\pi\sqrt{-1}\,(a_1\theta_1+\cdots+a_{s+1}\theta_{r+1})}\in S^1,
$$
where $S^1=\{z\in {\mathbb C}\vert\,\vert z\vert=1\}$, and 
$a_1,\,\,\ldots,\,a_{r+1}\in {\mathbb Z}$. 
The action of $T$ on $K\times S^1=SU(r+1)\times S^1$ 
by 
$$
(x,e^{2\pi\sqrt{-1}\,\theta})\cdot a=
(xa,\,\lambda_{\mathcal I}(a^{-1})\,e^{2\pi\sqrt{-1}\,\theta}),\qquad a\in T.
$$
The orbit space 
\begin{align*}
P={\mathcal S}_{\lambda}&=
SU(s+1)\times S^1/\!\!\sim\,\,\,\\
&=\{(x,e^{2\pi\sqrt{-1}\,\theta})\vert 
x\in SU(r+1),\,\,\theta\in {\mathbb R}\}/\!\!\sim
\end{align*}
whose equivalence relation is given by 
$(x',e^{2\pi\sqrt{-1}\,\theta'})\sim 
(x,e^{2\pi\sqrt{-1}\,\theta})$ 
is equivalent to that: 
  $x'=xt$ and $e^{2\pi\sqrt{-1}\,\theta'}=e^{2\pi\sqrt{-1}\,\theta}\lambda_{\mathcal I}(t^{-1})$. 
  We denote the equivalence class 
  including $(x,e^{2\pi\sqrt{-1}\,\theta})$ 
  by 
  $[x,e^{2\pi\sqrt{-1}\,\theta}]$. 
  Then, we have the principal $S^1$-bundle 
  $P={\mathcal S}_{\lambda}$ over 
  $K/T$ associated to 
  $\lambda_{\mathcal I}$, which is the space of all $T$-orbits through $(x,e^{2\pi\sqrt{-1}\,\theta})$, 
  $x\in SU(r+1),\,\,\theta\in {\mathbb R}$, namely, 
  $$
  P={\mathcal S}_{\lambda}=
  \{[x,e^{2\pi\sqrt{-1}\,\theta}]\vert\,x\in SU(r+1),\,\,\theta\in {\mathbb R}\}.
  $$
  \medskip\par
  {\bf Example 2.} \quad 
  In particular, let us consider the case $r=1$. Let 
  $$K=SU(2)\supset T
  =\bigg\{
  \begin{bmatrix}
  e^{2\pi\sqrt{-1}\,\theta}&0\\
  0&e^{-2\pi\sqrt{-1}\,\theta}
  \end{bmatrix}
  \big\vert\,\,\theta\in {\mathbb R}
  \bigg\},
  $$ 
  $\dim(K/T)=2$ and $\dim P=3$. 
  For $a_1,\,\,a_2\in {\mathbb Z}$, and $\ell=a_1-a_2$, 
  let 
  $$
  \lambda_{{\mathcal I}}:\,\,
  T\ni 
  \begin{bmatrix}
  e^{2\pi\sqrt{-1}\,\theta}&0\\
  0&e^{-2\pi\sqrt{-1}\,\theta}
  \end{bmatrix}
  \mapsto 
  e^{2\pi\sqrt{-1}\,(a_1-a_2)\theta}
  =e^{2\pi\sqrt{-1}\,\ell\,\theta}
  \in S^1
  $$
  and $T$ acts on $SU(2)\times S^1$ by 
  $$
  (x,e^{2\pi\sqrt{-1}\,\xi})\,\cdot\, a
  :=(xa,
  e^{2\pi\sqrt{-1}\,\ell\,\theta}\,
  e^{2\pi\sqrt{-1}\,\xi}
  ),
  $$
  for 
  $
  a=\begin{bmatrix}
  e^{2\pi\sqrt{-1}\,\theta}&0\\
  0&e^{-2\pi\sqrt{-1}\,\theta}
  \end{bmatrix}
  \in T,\,\, x\in SU(2),\,\,\xi\in {\mathbb R}
  $. 
  Then, $P$ is diffeomorphic with $S^3$, and 
  $M=K/T$ is diffeomorphic with $P^1\!({\mathbb C})$, 
  and we have 
  $\pi:\,\,P={\mathcal S}_{\lambda_{\mathcal I}}\rightarrow M=K/T=SU(2)/S^1=P^1\!({\mathbb C})$. 
  Let 
  \begin{align*}
  {\mathfrak k}&={\mathfrak su}(2)=\{X\in {\mathfrak gl}(2,{\mathbb C})\vert\,{}^{\mathrm t}\overline{X}+X=0,\,\,\mbox{\rm Tr}(X)=0\},\\
  {\mathfrak t}&={\mathfrak g}({\mathfrak u}(1)\times 
  {\mathfrak u}(1))
  =\left\{
  \begin{pmatrix}\sqrt{-1}\,\theta&0\\0&-\sqrt{-1}\,\theta
  \end{pmatrix}
  \vert\,\,\theta\in {\mathbb R}
  \right\},\\
  {\mathfrak m}&=
  \left\{\begin{pmatrix}
  0&-\overline{z}\\
  z&0
  \end{pmatrix}
  \vert\,\,z\in {\mathbb C}
  \right\},
  \end{align*}
  respectively. 
  Let $\langle\,\,\cdot\,,\,\cdot\,\,\rangle$ be 
  the inner product on $\frak k$ defined by 
  $$
  \langle X,Y\rangle:=-\frac12 \mbox{Tr}(X\,Y),
  \qquad X,\,\,Y\in {\frak k}. 
  $$
  Then, for $X=\begin{pmatrix}0&-\overline{z}\\z&0
  \end{pmatrix},\,\,
  Y=\begin{pmatrix}0&-\overline{w}\\w&0
  \end{pmatrix}\in {\frak m}$,
  $$
  \langle X,Y\rangle=
  x\xi+y\eta,\qquad 
  z=x+\sqrt{-1}\,y,\,\,w=\xi+\sqrt{-1}\,\eta,\,\,\,\,x,\,\,y,\,\,\xi,\,\,\eta\in {\mathbb R}, 
  $$
  and $h$, the $G$-invariant Riemannian metric 
  on $M=K/T= P^1\!({\mathbb C})$ 
  in such a way that 
  $$
  h_o(X_o,Y_o)=\langle X,Y\rangle,\qquad X,\,\,Y\in {\frak m},
  $$
  where $o=\{T\}\in M=K/T$. 
  Let $\{H_1,\,\,X_1,\,\,X_2\}$ be an orthonormal basis 
  of $\frak k$ with respect to $\langle\,\cdot\,,\,\cdot\,\rangle$ where 
  $$
  H_1=\begin{pmatrix}\sqrt{-1}&0\\0&-\sqrt{-1}
  \end{pmatrix},
  \quad 
  X_1=\begin{pmatrix}0&\sqrt{-1}\\\sqrt{-1}&0
  \end{pmatrix},
  \quad 
  X_2=\begin{pmatrix}0&-1\\
  1&0
  \end{pmatrix}
  $$
  satisfying that 
  $$
  [H_1,X_1]=2X_2,\quad 
  [X_2,H_1]=2X_1,\quad 
  [X_1,X_2]=2H_1.
  $$
  \par
  In our case, taking
$$
SU(2)\ni k\exp(sX_1+tX_2)\exp(uH_1)\mapsto 
(s,t,u)\in {\mathbb R}^3, 
$$
as a local coordinate around $k\in SU(2)$, and 
let us write a locally defined orthonormal frame field 
$\{e_i\}_{i=1}^3$ 
on $SU(2)$ around the identity $e$ in $SU(2)$ by 
$$
e_1=a\,\frac{\partial}{\partial s}+b\,\frac{\partial}{\partial t}, \,\,
e_2=c\,\frac{\partial}{\partial s}+d\,\frac{\partial}{\partial t}, \,\,
e_3=e^{C\ell(\ell-1)u\,(As+Bt)}\,\frac{\partial}{\partial u},
$$
where $a,\,\,b,\,\,c,\,\,d,\,\,A,\,\,B,\,\,C$ are real constants. 
\par
   For $X=\tau(\pi)=-(\kappa_1\widetilde{\epsilon}_1+\kappa_2\widetilde{\epsilon}_2)$, and 
   $\{e_i\}_{i=1}^3$ an orthonormal frame field on $P$ such that 
  the vertical subspace ${\mathcal V}_p={\mathbb R}e_{3\,p}$ and the horizontal subspace ${\mathcal H}_p={\mathbb R}\,e_{1\,p}\oplus {\mathbb R}\,e_{2\,p}$ 
  of $T_pP$ $(p\in P)$ satisfies 
  $$
  [e_i,\,e_3]=\kappa_i\,e_3\quad (i=1,2)
  $$
with $\kappa_i\in {C}^{\infty}(P)\,\, (i=1,2)$, where 
$\kappa_1=C\ell(\ell-1)u(aA+bB),\,\,\kappa_2=C\ell(\ell-1)u(cA+dB)$. 
It holds that 
\begin{equation}{\mathrm div}(X)=e_1\kappa_1+e_2\kappa_2\equiv 0.
\end{equation} 
Furthermore, we obtain 
\begin{align}
X=\tau(\pi)&=-(\kappa_1\widetilde{\epsilon}_1+\kappa_2\widetilde{\epsilon}_2)\nonumber\\
&=-C\ell(\ell-1)u\{(aA+bB)\widetilde{\epsilon}_1+(cA+dB)\widetilde{\epsilon}_2\}.
\end{align}
Therefore, if $\ell=0$ or $\ell=1$, 
$$
X=\tau(\pi)=0,
$$
namely, 
$\pi:\,\,P=S_{\lambda_{\mathcal I}}\rightarrow M=K/T=P^1\!({\mathbb C})$ is the direct product if $\ell=0$, and it is  
the standard Hopf fiberring is harmonic if $\ell=1$. 
\par 
If $\ell=2,3, \cdots$, 
our $X=\tau(\pi)\not\equiv0$ satisfies that 
$\overline{\Delta}^hX=cX$ with $\nabla^h_XX=0$ which is 
equivalent to 
$$
\overline{\Delta}^hX+\nabla^h_XX-{\mathrm Ric}^h(X)=0, 
$$
which is equivalent to that 
\begin{equation}
\overline{\Delta}^hX=cX,\quad \nabla^h_XX=0, 
\end{equation}
and $\pi:\,\,P={\mathcal S}_{\lambda_{\mathcal I}}\rightarrow M=K/T={\mathbb C}^1\!P$ is biharmonic, however it is not harmonic. 
Notice that $(M,h)=({\mathbb C}^1\!P,h)$ 
satisfies that 
${\mathrm Ric}^h=\frac12\,{\mathrm Id}$ with 
$c=\frac12$ and $\lambda_1(M,h)=1$ (\cite{TK}, p. 213,  and  \cite{U86}, p. 67, Type ${\mathrm A\,III}$ in Table A2 and also p. 70).  
\vskip0.3cm\par
Therefore, we can summarize: 
\begin{th}
For $\ell=1,2,\ldots$, let $$
  \lambda_{{\mathcal I}}:\,\,
  T\ni 
  \begin{bmatrix}
  e^{2\pi\sqrt{-1}\,\theta}&0\\
  0&e^{-2\pi\sqrt{-1}\,\theta}
  \end{bmatrix}
  \mapsto 
  e^{2\pi\sqrt{-1}\,\ell\,\theta}
  \in S^1
  $$
be a homomorphism of $T$ into $S^1$, and let 
$\pi:\,\,P={\mathcal S}_{\lambda_{\mathcal I}}\rightarrow M=K/T=SU(2)/S^1=P^1\!({\mathbb C})$ be the 
principal $S^1$-bundle over $K/T$ associated to 
$\lambda_{\mathcal I}$.  
Then, for every $\ell=2,3,\ldots$, the projection 
$\pi:\,(P,g)\rightarrow (M,h)$ is biharmonic but not harmonic. 
\end{th}


\begin{thebibliography}{99}
\bibitem{AM} K. Akutagawa and S. Maeta, 
\textit{Complete biharmonic submanifolds in the Euclidean spaces}, 
Geometriae Dedicata, \textbf{164} (2013), 351--355. 
\bibitem{AO} M.A. Akyol and Y.L. Ou, 
\textit{Biharmonic Riemannian submersions}, 2018, arXiv: 1805.04754v1. 
\bibitem{BK} P. Paird and D. Kamissoko, 
\textit{On constructing biharmonic maps and metrics},  
Ann. Global Anal. Geom. \textbf{23} (2003), 65--75. 
\bibitem{BW} P. Baird and J.C. Wood, {\em Harmonic morphisms between Riemannian manifolds}, London Math. Soc. Monographs, Oxford, 2003.  
\bibitem{BMO1}
A. Balmus, S. Montaldo and C. Oniciuc, 
\textit{Classification results for biharmonic submanifolds in spheres}, Israel J. Math., 
\textbf{168} (2008), 201--220. 
\bibitem{BMO2} 
A. Balmus, S. Montaldo and C. Oniciuc, 
\textit{Biharmonic hypersurfaces in $4$-dimensional space forms}, Math. Nachr., \textbf{283} (2010), 1696--1705. 
\bibitem{B} R. Bott, 
\textit{Homogeneous vector bundles}, 
Ann. Math., \textbf{66}, 203--248. 
\bibitem{BG} C. Boyer and K. Galicki, 
\textit{Sasakian Geometry}, 
Oxford Sci. Publ., 2008. 
\bibitem{CMP} R. Caddeo, S. Montaldo, P. Piu, 
\textit{On biharmonic maps}, Contemp. Math., \textbf{288} (2001), 286--290. 
\bibitem{CLU} 
I. Castro, H.Z. Li and F. Urbano, 
\textit{Hamiltonian-minimal Lagrangian submanifolds in complex space forms}, 
Pacific J. Math., \textbf{227} (2006), 43--63. 
\bibitem{C} B.Y. Chen, \textit{Some open problems and conjectures on submanifolds of finite type}, 
Soochow J. Math., \textbf{17} (1991), 169--188. 
\bibitem{D} 
F. Defever, \textit{Hypersurfaces in ${\mathbb E}^4$ with harmonic mean curvature vetor}, Math. Nachr., {\bf 196} (1998), 61--69. 
\bibitem{EL1}
J. Eells and L. Lemaire, \textit{Selected Topics in Harmonic Maps}, 
CBMS, Regional Conference Series in Math., Amer. Math. Soc., {\bf 50}, 1983.  
\bibitem{FO} 
D. Fetcu and C. Oniciuc, 
\textit{Biharmonic integral ${\mathcal C}$-parallel submanifolds in 7-dimensional Sasakian space forms}, 
Tohoku Math. J., {\bf 64} (2012), 195--222.  
\bibitem{HV} 
T. Hasanis and T. Vlachos 
\textit{Hypersurfaces in ${\mathbb E}^4$ with harmonic mean curvaturer vector field}, 
Math. Nachr., {\bf 172} (1995), 145--169. 
\bibitem{IIU2} 
T. Ichiyama, J. Inoguchi, H. Urakawa, \textit{Biharmonic maps and bi-Yang-Mills fields}, 
Note di Mat., {\bf 28}, (2009), 233--275. 
\bibitem{IIU} T. Ichiyama, J. Inoguchi, H. Urakawa, \textit{Classifications and isolation phenomena  of biharmonic maps and bi-Yang-Mills fields}, 
Note di Mat., {\bf 30}, (2010), 15--48. 
\bibitem{I1} 
J. Inoguchi, \textit{Submanifolds with harmonic mean curvature vector filed in contact 3-manifolds}, Colloq. Math., \textbf{100} (2004), 163--179. 
\bibitem{Ir} H. Iriyeh, 
\textit{Hamiltonian minimal Lagrangian cones in ${\mathbb C}^m$}, 
Tokyo J. Math., \textbf{28} (2005), 91--107. 
\bibitem{II} S. Ishihara and S. Ishikawa, \textit{Notes on relatively harmonic immersions}, Hokkaido Math. J., \textbf{4} (1975), 234--246. 
\bibitem{J} G.Y. Jiang, \textit{2-harmonic maps and their first and second variational formula}, Chinese Ann. Math., \textbf{7A} (1986), 
388--402;  Note di Mat., {\bf 28} (2009), 209--232.
\bibitem{K1} T. Kajigaya, {\em Second variation formula and the stability of Legendrian minimal submanifolds in Sasakian manifolds}, 
Tohoku Math. J., \textbf{65} (2013), 523--543. 
\bibitem{K} S. Kobayashi, \textit{Transformation Groups in Differential Geometry}, Springer, 1972. 
\bibitem{LO1} E. Loubeau, C. Oniciuc, \textit{The index of biharmonic maps in spheres}, 
Compositio Math., \textbf{141} (2005), 729--745. 
\bibitem{LO2} E. Loubeau and C. Oniciuc, \textit{On the biharmonic and harmonic indices of the Hopf map}, Trans. Amer. Math. Soc., {\bf 359} (2007), 5239--5256. 
\bibitem{LOu}
E. Loubeau and Y-L. Ou, 
\textit{
Biharmonic maps and morphisms from conformal mappings}, 
Tohoku Math. J., {\bf 62} (2010), 55--73. 
\bibitem{MU} S. Maeta and U. Urakawa, 
\textit{Biharmonic Lagrangian submanifolds in K\"ahler manifolds}, 
Glasgow Math. J. , {\bf 55} (2013), 465--480. 
\bibitem{MO1} S. Montaldo, C. Oniciuc, \textit{A short survey on biharmonic maps between Riemannian manifolds}, Rev. Un. Mat. Argentina \textbf{47} (2006), 1--22. 
\bibitem{Ng} 
Y. Nagatomo, 
\textit{Harmonic maps into Grassmannians and a generalization of do Carmo-Wallach theorem}, Proc. the 16th OCU Intern. Academic Symp. 2008, 
OCAMI Studies, \textbf{3} (2008), 41--52.  
\bibitem{NU1} 
N. Nakauchi and H. Urakawa, 
\textit{Biharmonic hypersurfaces in a Riemannian manifold with non-positive Ricci curvature}, 
Ann. Global Anal. Geom., \textbf{40} (2011), 125--131. 
\bibitem{NU2} 
N. Nakauchi and H. Urakawa, 
\textit{Biharmonic submanifolds in a Riemannian 
manifold with non-positive curvature}, Results in Math.,{\bf 63} (2013), 467--474. 
\bibitem{NUG} 
N. Nakauchi, H. Urakawa and S. Gudmundsson, 
\textit{Biharmonic maps into a Riemannian manifold of non-positive curvature}, 
Geom. Dedicata, \textbf{169}. (2014), 263--272.   
\bibitem{O} 
C. Oniciuc, \textit{Biharmonic maps between Riemannian manifolds}, 
Ann. Stiint Univ. A${\ell}$. I. Cuza Iasi, Mat. (N.S.), {\bf 48} No. 2, (2002), 
237--248. 
\bibitem{ON} B. O'Neill, \textit{The fundamental equation of a submersion}, Michigan Math. J., {\bf 13} (1966), 459--469.  
\bibitem{OSU} S. Ohno, T. Sakai and H. Urakawa, 
{\em Biharmonic homogeneous hypersurfaces in compact symmetric spaces}, Differ. Geom. Appl., {\bf 43} (2015), 155--179. 
\bibitem{OSU2} S. Ohno, T. Sakai and 
H. Urakawa, 
{\em Rigidity of transversally biharmonic maps between foliated Riemannian manifolds}, to appear in Hokkaido Math. J. 
\bibitem{OT} 
Ye-Lin Ou and Liang Tang, 
{\em The generalized Chen's conjecture on biharmonic submanifolds is false}, 
arXiv: 1006.1838v1.  
\bibitem{OT2} 
Ye-Lin Ou and Liang Tang, \textit{On the generalized Chen's conjecture on biharmonic submanifolds}, Michigan Math. J., {\bf 61} (2012), 531--542. 
\bibitem{S1} T. Sasahara, \textit{Legendre surfaces in Sasakian space forms whose mean curvature vectors are eigenvectors}, Publ. Math. 
Debrecen, \textbf{67} (2005), 285--303. 
\bibitem{S2} T. Sasahara, \textit{Stability of biharmonic Legendrian submanifolds in Sasakian space forms}, 
Canad. Math. Bull. \textbf{51} (2008), 448--459. 
\bibitem{S3} T. Sasahara, \textit{A class of biminimal Legendrian submanifolds in Sasaki space forms}, a preprint, 2013, to appear in Math. Nach. 
\bibitem{T}
T. Takahashi, \textit{Minimal immersions of Riemannian manifoplds}, 
J. Math. Soc. Japan, \textbf{18} (1966), 380--385. 
\bibitem{TK} M. Takeuchi and S. Kobayashi, 
\textit{Minimal imbeddings of $\mathbb R$-space}, 
J. Differ. Geom., \textbf{2} (1968), 203--213. 
\bibitem{U86} H. Urakawa, \textit{The first eigenvalue of the Laplacian for a positively curved homogeneous Riemannian manifold}, Compositio Math., \textbf{59} (1986), 57--71.  
\bibitem{U} H. Urakawa, \textit{CR rigidity of pseudo harmonic maps and pseudo biharmonic maps}, Hokkaido Math. J., \textbf {46} (2017), 141--187.  
\bibitem{U1} H. Urakawa, \textit{Harmonic maps and biharmonic maps on principal bundles and warped products}, to appear in J. Korean Math. Soc., 2018. 
\bibitem{U2} H. Urakawa, \textit{Calculus of Variations and Harmonic Maps}, Vol. 132, Amer. Math. Soc., 1990.  
\bibitem{WO} 
Z-P Wang and Y-L Ou, \textit{Biharmonic Riemannian submersions from 3-manifolds}, Math. Z., \textbf{269} (2011), 917--925.  
\end{thebibliography}
\end{document}